\documentclass[10pt,twoside]{article}
\usepackage{times, fancyhdr}
\usepackage{cite}
\usepackage[ ]{graphicx}
\usepackage{amsmath}
\usepackage{amssymb}
\usepackage{amsfonts}
\usepackage{latexsym,amsmath,amsfonts,amssymb,amsthm}

\usepackage{amsthm}
\usepackage{indentfirst}
\usepackage{sectsty}
\sectionfont{\large}
\subsectionfont{\normalsize\bf}

\usepackage{float}

\numberwithin{equation}{section}

\newtheorem{theorem}{Theorem}[section]

\newtheorem{cor}[theorem]{Corollary}

\begin{document}

%For IJMM Editorial Office Use Only.
%%%%%%%%%%%%%%%%%%%%%%%%%%%%%%%%%%%

\title{\vspace{-0.844in}\parbox{\linewidth}{\footnotesize\noindent
{Int. J. Mod. Math. 2(2007), no.\,2, 247--251.}
\hspace{\stretch{1}}}
%\copyright 2007 Dixie W publishing Corporation, U. S. A.}
%\vspace{\bigskipamount} \\
%%%%%%%%%%%%%%%%%%%%%%%%%%%%%%%%%%%
{\large An Extension of a Curious Binomial Identity}}
\date{}

\author{\normalsize Zhi-Wei Sun and Ke-Jian Wu}
%\\ {\footnotesize Received April 10, 2007; Revised June 10, 2007} }
\maketitle

\thispagestyle{empty}

\renewcommand{\theequation}{\thesection.\arabic{equation}}
\makeatletter
\@addtoreset{equation}{section}

\date{}{}
\def\colon{{:}\;}
\def\Z{\Bbb Z}
\def\zp{\Z^+}
\def\Q{\Bbb Q}
\def\N{\Bbb N}
\def\R{\Bbb R}
\def\C{\Bbb C}
\def\al{\alpha}
\def\l{\left}
\def\r{\right}
\def\bg{\bigg}
\def\({\bg(}
\def\){\bg)}
\def\[{\bg\lfloor}
\def\]{\bg\rfloor}
\def\t{\mbox}
\def\f{\frac}
\def\si{\sim_2}
\def\pr{\prec_2}
\def\pe{\preceq_2}
\def\B{B_{k+1}}
\def\E{E_k}
\def\P{\Z_q}
\def\q{\ (\t{\rm mod}\ q)}
\def\qc{q\t{-adic}}
\def\se {\subseteq}
\def\sp {\supseteq}
\def\sm{\setminus}
\def\Ar{\Longrightarrow}
\def\La{\Longleftrightarrow}
\def\bi{\binom}
\def\eq{\equiv}
\def\cs{\ldots}
\def\ls{\leqslant}
\def\gs{\geqslant}
\def\mo{\t{\rm mod}}
\def\Tor{\t{\rm Tor}}
\def\sign{\t{\rm sign}}
\def\o{\t{\rm ord}}
\def\ch{\t{\rm ch}}
\def\per{\t{\rm per}}
\def\ve{\varepsilon}
\def\da{\delta}
\def\Da{\Delta}
\def\la{\lambda}
\def\ta{\theta}
\def\ord{\t{\rm ord}}
\def\Tor{\t{\rm Tor}}
\def\kl{k\varphi(q^{\al})+r}
\def\mx{\langle-x\rangle_q}
\def\Mr{\sum\Sb 0\ls j<q\\m\mid j-r-qs\endSb(x+j)^k}
\def\qed{}
\def\Proof{\noindent{\it Proof}}
\def\Def{\medskip\noindent{\it Definition}\ }
\def\Remark{\medskip\noindent{\it  Remark}}
\def\Ack{\medskip\noindent {\bf Acknowledgments}}

\begin{abstract}
\noindent In 2002 Z. W. Sun published a curious identity involving
binomial coefficients. In this paper we obtain the following
generalization of the identity:
\begin{align*}
&\l(x+(m+1)z\r)\sum^m_{n=0}(-1)^n\bi {x+y+nz}{m-n}\bi {y+n(z+1)}n\\
=&z\sum_{0\leqslant l\leqslant n\leqslant m}(-1)^n\bi nl\bi
{x+l}{m-n}(1+z)^{n+l}(1-z)^{n-l}+(x-m)\bi xm.
\end{align*}

\smallskip
\smallskip
\smallskip
\noindent {\bf Keywords:} Binomial coefficient, combinatorial
identity.

\smallskip
\noindent {\bf 2000 Mathematics Subject Classification:} 05A19, 11B65.
\end{abstract}

\section{Introduction}

In 2002 Z. W. Sun \cite{S} established the following combinatorial
identity:
\begin{equation}\label{eq-1.1}
\aligned&(x+m+1)\sum ^m_{n=0}(-1)^n\bi {x+y+n}{m-n}\bi {y+2n}n
\\&\ =\sum^m_{n=0}\bi {x+n}{m-n}(-4)^n+(x-m)\bi xm,
\endaligned
\end{equation}
where $m\in\N=\{0,1,2,\cs\}.$ Later A. Panholzer and H. Prodinger
\cite{PP} gave a new proof using generating functions, D. Merlini
and R. Sprugnoli \cite{MS} obtained another proof by means of
Riordan arrays, S. B. Ekhad and M. Mohammed \cite{EM} presented a
WZ proof of the identity, and W. Chu and L.V.D. Claudio \cite{CC}
re-proved the identity by using Jensen's formula.

In this paper we aim to extend the curious identity as follows.

\begin{theorem} For any $m=0,1,2,\cs$, we have
\begin{equation}\label{eq-1.2}\aligned
&\l(x+(m+1)z\r)\sum^m_{n=0}(-1)^n\bi {x+y+nz}{m-n}\bi {y+n(z+1)}n\\
=&z\sum_{0\leqslant l\leqslant n\leqslant m}(-1)^n\bi nl\bi
{x+l}{m-n}(1+z)^{n+l}(1-z)^{n-l}+(x-m)\bi xm;
\endaligned
\end{equation}
equivalently,
\begin{equation}\label{eq-1.3}\aligned
&\l(x+(m+1)z+1\r)\sum^m_{n=0}(-1)^n\bi {x+y+nz}{m-n}\bi {y+n(z+1)}n\\
=&(z+1)\sum_{0\leqslant l\leqslant n\leqslant m}(-1)^n\bi nl\bi
{x+l+1}{m-n}(1+z)^{n+l}(1-z)^{n-l}+(x-m)\bi xm.
\endaligned
\end{equation}
\end{theorem}

\Remark\ 1.1. Soon after the initial version of this paper was
posted as a preprint ({\tt arXiv:math.CO/0401057}), D. Callan [C]
found a nice combinatorial interpretation of the identity
(\ref{eq-1.1}) (which was called ``Sun's identity" by him) and
also a slightly more complicated combinatorial proof of our
generalization (\ref{eq-1.2}).

\medskip
Clearly, (\ref{eq-1.2}) in the case $z=1$ gives Sun's identity (\ref{eq-1.1}), and
(\ref{eq-1.3}) in the case $z=1$ yields the following equivalent form of (\ref{eq-1.1}).

\begin{cor} Let $m\in\N$. Then
\begin{equation}\label{eq-1.4}
\aligned&(x+m+2)\sum ^m_{n=0}(-1)^n\bi {x+y+n}{m-n}\bi {y+2n}n
\\&\ \ =2\sum^m_{n=0}\bi {x+n+1}{m-n}(-4)^n+(x-m)\bi xm.
\endaligned
\end{equation}
\end{cor}

\Remark\ 1.2. (\ref{eq-1.4}) in the special case $x-m\in\N$ and $y=1$, was ever conjectured by Z. H. Sun.

\begin{cor} For any $m\in\N$ we have
$$\sum_{0\leqslant l\leqslant n\leqslant m}(-1)^n\bi nl
\bi {l+(m+1)z}{m-n}(1+z)^{n-l}(1-z)^{n+l}
=(m+1)\bi {(m+1)z-1}m.$$
\end{cor}

\Proof. Just take $x=-(m+1)z$ in (\ref{eq-1.2}) and then replace $z$ by $-z$. \qed

\section{Proof of Theorem 1.1}

The starting point of our proof of Theorem 1.1 is the following known identity:
\begin{equation}\label{eq-2.1}
\sum_{n=0}^\infty\bi{\al+n\beta}n\l(\f{x-1}{x^\beta}\r)^n=\f{x^{\al+1}}{(1-\beta)x+\beta}.
\end{equation}
It appeared as (9) of H. W. Gould \cite{G}, and dates back to an identity
of Lambert (cf. (E.3.1) of \cite{AAR}).
Both (\ref{eq-2.1}) and Lambert's identity can be proved by Lagrange's inversion formula
(see pp. 631--632 of \cite{AAR}). In 2005 V. J. W. Guo and J. Zeng \cite{GZ} applied (\ref{eq-2.1})
to deduce some combinatorial identities originally motivated by the enumeration of convex polyominoes.

Let $\C$ be the complex field. For a formal power series
$f(t)\in\C[\![t]\!]$, the coefficient of
$t^n$ in $f(t)$ will be denoted by $[t^n]f(t)$.

\medskip

\noindent{\it Proof of Theorem 1.1}. In the case $m=0$, both (\ref{eq-1.2}) and
(\ref{eq-1.3}) are trivial. Below we assume that $m$ is a positive integer.

 Putting $\al=y$, $\beta=z+1$ and $x=1/(1+t)$ in (\ref{eq-2.1}), we find that
$$\sum^{\infty}_{n=0}\bi
{y+n(z+1)}n\l(-t(1+t)^z\r)^n=\f{(1+t)^{-y}}{1+t(z+1)}.$$
Thus
\begin{align*}
[t^m]\f {(1+t)^x}{1+t(z+1)}=&[t^m]\sum^{\infty}_{n=0}\bi {y+n(z+1)}n(-t)^n(1+t)^{nz+x+y}\\
=&\sum^m_{n=0}(-1)^n\bi {y+n(z+1)}n\bi {x+y+nz}{m-n}.
\end{align*}
(As pointed out by one of the referees, this identity
can be reproved by a mixed use of Lagrange's inversion formula and the Riordan array method.)
On the other hand,
\begin{align*}
&[t^m]\f {(1+t)^x}{(1+t(z+1))^2}=[t^m]\f {(1+t)^x}{1+t(z+1)(t(z+1)+2)}\\
=&[t^m](1+t)^x\sum^{\infty}_{n=0}\l(-t(z+1)\l(t(z+1)+2\r)\r)^n\\
=&\sum^m_{n=0}(-1)^n(z+1)^n[t^{m-n}](1+t)^x\l(t(z+1)+2\r)^n\\
=&\sum^m_{n=0}(-1)^n(z+1)^n[t^{m-n}](1+t)^x\l((z+1)(1+t)+1-z\r)^n\\
=&\sum^m_{n=0}(-1)^n(z+1)^n[t^{m-n}](1+t)^x\sum^n_{l=0}\bi nl(z+1)^l(1+t)^l(1-z)^{n-l}\\
=&\sum^m_{n=0}(-1)^n(z+1)^n\sum^n_{l=0}\bi nl\bi {x+l}{m-n}(1+z)^l(1-z)^{n-l}\\
=&\sum_{0\leqslant l\leqslant n\leqslant m}(-1)^n\bi nl\bi
{x+l}{m-n}(1+z)^{n+l}(1-z)^{n-l}.
\end{align*}

Since
\begin{align*}
[t^m]\f {(1+t)^x}{1+t(z+1)}=&[t^m]\f {(1+t)^x}{(1+t(z+1))^2}((1+t)(z+1)-z)\\
=&(z+1)[t^m]\f {(1+t)^{x+1}}{(1+t(z+1))^2}-z[t^m]\f
{(1+t)^x}{(1+t(z+1))^2},
\end{align*}
by the above we have
\begin{align*}
&\sum^m_{n=0}(-1)^n\bi {x+y+nz}{m-n}\bi {y+n(z+1)}n\\
=&(z+1)\sum_{0\leqslant l\leqslant n\leqslant m}
(-1)^n\bi nl\bi {x+1+l}{m-n}(1+z)^{n+l}(1-z)^{n-l}\\
&-z\sum_{0\leqslant l\leqslant n\leqslant m}(-1)^n\bi nl\bi
{x+l}{m-n}(1+z)^{n+l}(1-z)^{n-l}.
\end{align*}
From this we immediately see that (\ref{eq-1.2}) and (\ref{eq-1.3}) are equivalent.

Observe that
\begin{align*}
m[t^m]\f {(1+t)^x}{1+t(z+1)}=&[t^{m-1}]\f {\partial}{\partial t}
\(\f {(1+t)^x}{1+t(z+1)}\)\\
=&[t^{m-1}]\(\f {-(z+1)(1+t)^x}{(1+t(z+1))^2}+\f{x(1+t)^{x-1}}{1+t(z+1)}\)\\
=&[t^m]\(-\f {t(z+1)(1+t)^x}{(1+t(z+1))^2}+\f{xt(1+t)^{x-1}}{1+t(z+1)}\)\\
=&[t^m]\(\f {(1+t)^x}{(1+t(z+1))^2}+\f
{xt(1+t)^{x-1}-(1+t)^x}{1+t(z+1)}\)
\end{align*}
and hence
$$(m+1)[t^m]\f {(1+t)^x}{1+t(z+1)}-[t^m]\f {(1+t)^x}{\big(1+t(z+1)\big)^2}
=[t^m]\f {xt(1+t)^{x-1}}{1+t(z+1)}.$$
It follows that
\begin{align*} &\big(x+(m+1)z\big)[t^m]\f {(1+t)^x}{1+t(z+1)}
-z[t^m]\f {(1+t)^x}{\big(1+t(z+1)\big)^2}
\\=&x[t^m]\f {(1+t)^x}{1+t(z+1)}+z[t^m]\f {xt(1+t)^{x-1}}{1+t(z+1)}
=x[t^m]\f {(1+t)^{x-1}(1+t+zt)}{1+t(z+1)}
\\=&x[t^m](1+t)^{x-1}=x\bi {x-1}m=(x-m)\bi xm.
\end{align*}
This, together with the previous arguments, yields the identity (\ref{eq-1.2}).

The proof of Theorem 1.1 is now complete. \qed

\Ack. The first author is supported by a Key
Program (Grant No. 10331020) of the National Natural Science Foundation in China.
Both authors thank the referees for helpful comments.

\vspace{5mm}

Zhi-Wei Sun: Department of Mathematics, Nanjing University, Nanjing 210093,
People's Republic of China

{\it E-mail address:} {\tt zwsun@nju.edu.cn}

\vspace{2mm}

Ke-Jian Wu: Department of Mathematics, Zhanjiang Normal College,
Zhanjiang 524048, People's Republic of China

{\it E-mail address:} {\tt kjwu328@yahoo.com.cn}

\end{document}